\documentclass[11pt]{article}
\usepackage{amsfonts}
\usepackage{amsmath,amsthm,color,epsfig}
\usepackage{amssymb}
\usepackage{amsbsy,mathtools}
\usepackage{enumerate,caption}
\title{Superalgebraically Convergent Smoothly-Windowed Lattice Sums
  for Doubly Periodic Green Functions in Three-Dimensional Space}

\author{Oscar P. Bruno\footnote{Applied and Computational Mathematics,
    Caltech, Pasadena, CA 91125. Email: obruno@caltech.edu},
  \, Stephen P. Shipman\footnote{Dept. of Mathematics, Louisiana State
    University, Baton Rouge, LA \ 70803. Email:
    shipman@math.lsu.edu},
  \, Catalin Turc\footnote{Dept. of Math.
    Sciences, New Jersey Inst. of Technology, Newark, NJ 07102. Email: catalin.c.turc@njit.edu},
    \, Stephanos
  Venakides\footnote{Dept. of Mathematics, Duke University, Durham, NC
    \ 27708. Email: ven@math.duke.edu} }
\newtheorem{Theorem}{Theorem}[section]

\newcommand{\xx}{{\hat x}}
\newcommand{\yy}{{\hat y}}
\newcommand{\zz}{{\hat z}}
\newcommand{\rr}{\boldsymbol{r}}
\newcommand{\xxi}{\boldsymbol{\xi}}

\newcommand{\bfx}{{\mathbf{x}}}

\newcommand{\veps}{\varepsilon}

\newcommand{\ZZ}{\mathbb{Z}}
\newcommand{\RR}{\mathbb{R}}

\newcommand{\vv}{\mathbf v}

\begin{document}
\date{}
\maketitle
\begin{abstract}
  This paper, Part~I in a two-part series, presents {\em (i)}~A simple and highly efficient algorithm for evaluation of quasi-periodic Green functions, as well as {\em (ii)}~An associated boundary-integral equation method for the numerical solution of problems of  scattering of waves by doubly periodic arrays of scatterers in three-dimensional space.  Except for certain ``Wood frequencies'' at  which the quasi-periodic Green function ceases to exist, the  proposed approach, which is based on use of smooth windowing  functions, gives rise to lattice sums which converge  superalgebraically fast---that is, faster than any power of the  number of terms used---in sharp contrast with the extremely slow  convergence exhibited by the corresponding sums in absence of smooth  windowing. (The Wood-frequency problem is treated in Part~II.\footnote{A preliminary version of Part II can be found in section~4 of arXiv:1307.1176v1.}) A  proof presented in this paper establishes rigorously the superalgebraic convergence of the windowed lattice sums.  A variety  of numerical results demonstrate the practical efficiency  of the proposed approach.

\medskip
\indent $\mathbf{Keywords}$: scattering, periodic Green function, lattice sum, smooth truncation, super-algebraic convergence, boundary-integral equations.
\end{abstract}

\newpage


\section{Introduction}

The numerical solution of problems of electromagnetic, acoustic and
elastic wave scattering by doubly periodic structures entails
significant difficulties.  Assuming harmonic temporal dependence with
frequency $\omega$, the scattered fields can be obtained by means of
numerical methods based on integral equations---provided a
corresponding numerical scheme is used to evaluate the classical
radiating quasi-periodic Green function $G^{qper}$ for the
three-dimensional scalar Helmholtz operator $H[u] = \Delta u + k^2 u$
($k = \omega/c$ where $c$ is the propagation speed). The
aforementioned difficulties arise, to a significant extent, from
corresponding challenges posed by the evaluation of the quasi-periodic
Green function.

The quasi-periodic Green function $G^{qper}$ can be constructed as an
infinite sum of free-space Green functions (Helmholtz monopoles) with
bi-periodically distributed monopole singularities.
Let $\vv_1$ and $\vv_2$ denote two independent vectors in $\RR^2$
that characterize the periodicity, and let $\vv_1^*$ and $\vv_2^*$ be the dual vectors, that is $\vv_i^*\cdot\vv_j=\delta_{ij}$.  The Bloch wavevector will be denoted by $\mathbf{k}=\alpha\vv_1^*+\beta\vv_2^*$, where $\alpha$ and $\beta$ are the Bloch wavenumbers.
With the notation
$\bfx= (x,y,z)\in \RR^3$ and $\tilde\bfx=(x,y)$ and
\begin{equation}\label{r}
 r_{mn}^2= \left| \tilde\bfx + m\vv_1 + n\vv_2 \right|^2 +z^2,
\end{equation}
the quasi-periodic Green function
can be expressed in the form
\begin{equation}\label{latticesum1}
G^{qper}(x,y,z) \,=\, \frac{1}{4\pi}\sum_{m,n\in\ZZ} \frac{e^{ikr_{mn}}}{r_{mn}}
  \cdot e^{-i\mathbf{k}\cdot( m \mathbf{v}_1+ n \mathbf{v}_2)}.
\end{equation}
Notice that $\mathbf{k}\cdot( m \mathbf{v}_1+ n \mathbf{v}_2) = \alpha m + \beta n$.
The function $G^{qper}(x,y,z)$ possesses the quasi-periodic property
\begin{equation}
  G^{qper}(\tilde\bfx + m\vv_1+n\vv_2,z)
  = G^{qper}(\tilde\bfx,z) e^{i(\alpha m + \beta n)}.
\end{equation}

The series expansion (\ref{latticesum1}) possesses notoriously poor convergence
properties.  Various methods to accelerate its convergence, notably
the Ewald method~\cite{Ewald,Capolino,Papanicolaou}, have been
proposed.  A survey in these regards is given in~\cite{Linton}, and a
comprehensive discussion of lattice summation techniques can be found
in~\cite{Encyclopedia}. A few remarks concerning the computational
costs associated with previous accelerated methods for evaluation of
the Green function~\eqref{latticesum1} are presented below in this
section.

In the approach proposed presently, the infinite
sum~(\ref{latticesum1}) is evaluated by multiplying its $(m,n)$-th
term by the value $\chi_a(\tilde r_{mn})$ of a slow-rise smooth
windowing function $\chi_a$ which, evaluated at the {\em cylindrical
  radius}
\begin{equation}\label{rtilde}
  \tilde r_{mn}=\left| m\vv_1 + n\vv_2 \right|,
\end{equation}
restricts the sum to values of $m$ and $n$ satisfying $0\leq \tilde
r_{mn}\leq a$. (Note that $\tilde r_{mn} = r_{mn}$ if and only if
$z=0$.)  The function $\chi_a = \chi_a (\tilde r)$ is obtained as a
scaled version of an infinitely smooth real valued function $\chi$
defined on the set of non-negative real numbers $\tilde r\geq 0$,
which equals zero for $\tilde r>a$ and equals $1$ for $0\leq \tilde
r<c\,a$, where $c<1$ is an adequately selected real number. (For the
numerical experiments presented in this paper the value $c=0.5$ was
used.) The function $\chi_a$ is then defined by
\begin{equation}\label{chi_def}
\chi_a(\tilde r)= \chi(\tilde r/a);
\end{equation}
clearly $\chi_a$ decreases from $1$ to $0$ in a slow manner: its
derivative tends to zero as $a\to\infty$ throughout the region of decrease
$c\cdot a\leq \tilde r\leq a$.

\begin{figure} 
\begin{center}
\scalebox{0.6}{\includegraphics{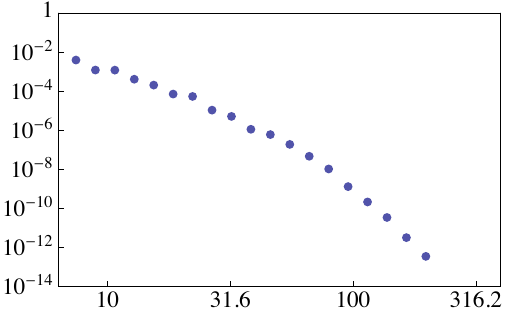}}
\scalebox{0.6}{\includegraphics{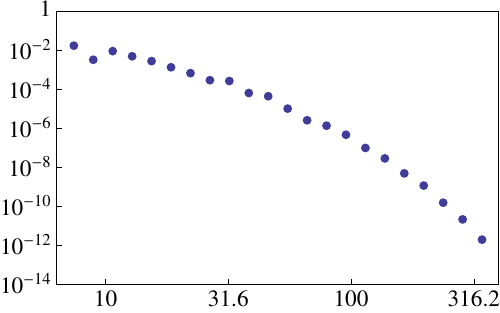}}
\scalebox{0.6}{\includegraphics{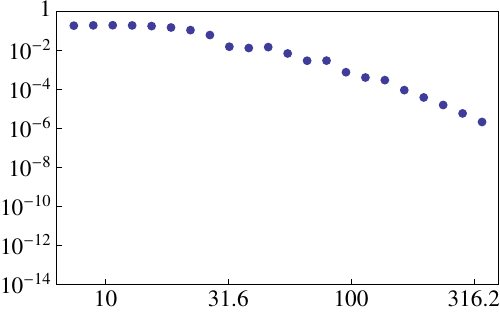}}

\scalebox{0.6}{\includegraphics{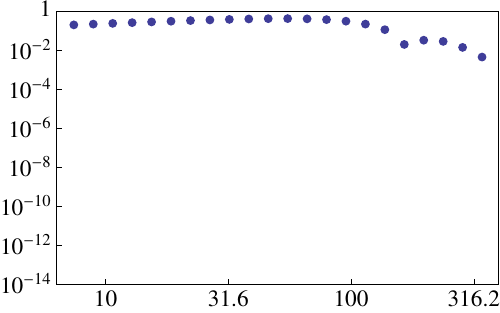}}
\scalebox{0.6}{\includegraphics{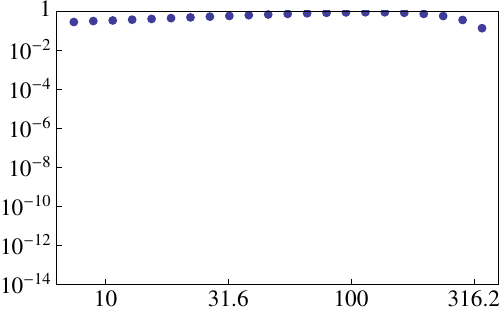}}
\scalebox{0.6}{\includegraphics{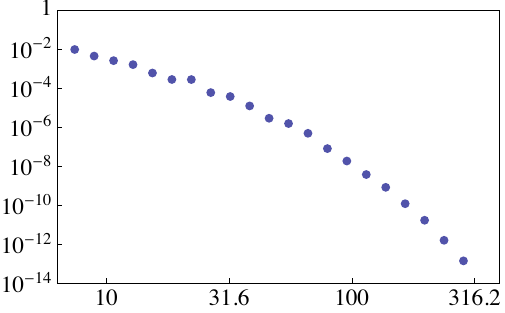}}
\end{center}
\caption{\small The error in the approximation of the quasi-periodic Green function by multiplying the lattice sum by a smooth truncation function $\chi((m+x)/a)\chi((n+y)/a)$, in which
$\chi(s)=\exp(2e^{1/(1-x)}/(x-2))$.
The plots show $\max_{(x,y,z)\in K}|G_{i+1}-G_i|$ as a function of $a_i$ on a $log$-$log$ scale, in which a truncated lattice sum $G_i$ is computed for $a=a_i=1.2^i$, $\hat{\mathbf x}=(\xx,\yy,\zz)=(0,0,1)$ and $(x,y,z)$ on a grid $K$ of evenly spaced points in $[0,0.6]\times[0,0.6]\times[0.6,1.4]$, excluding $\hat{\mathbf x}={\mathbf x}$.
The lattice vectors are $\mathbf{v}_1=(1,0,0)$ and $\mathbf{v}_2=(0,1,0)$, the Bloch wavevector is $(\kappa_1=0,\kappa_2=0)$, and the frequencies are $k=0.4,0.8,0.95$ (first row) and $k=0.99,2.24,2.5$ (second row).  Both $k=1.0$ and $k\approx2.23607$ are Wood frequencies, at which convergence is not available.}
\label{fig:GreenConvergenceNormal}
\end{figure}

\begin{figure}
\begin{center}
\scalebox{0.6}{\includegraphics{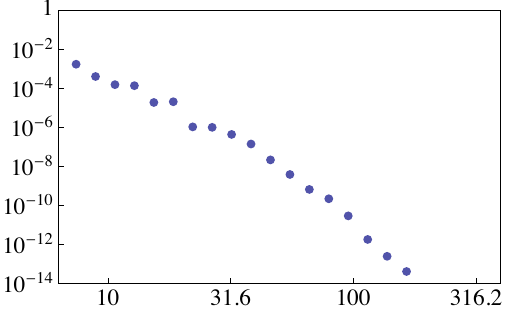}}
\scalebox{0.6}{\includegraphics{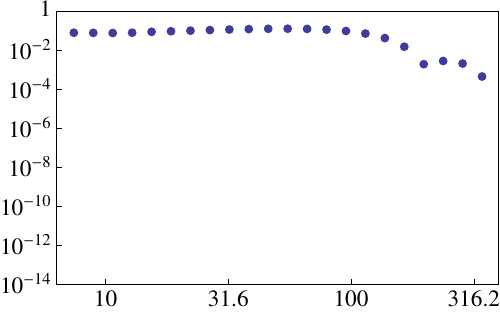}}
\scalebox{0.6}{\includegraphics{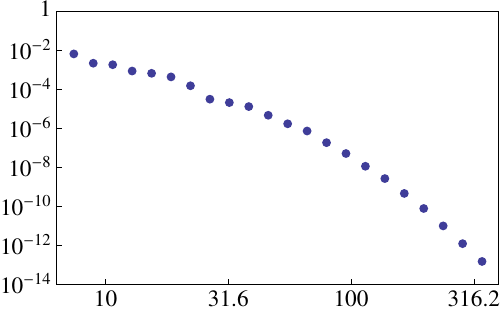}}
\end{center}
\caption{\small These plots are similar to those in Fig.~\ref{fig:GreenConvergenceNormal} except that the Bloch wavevector is $(\kappa_1=0.4,\kappa_2=-0.3)$, and the frequencies are $k=0.3,0.93,1.1$.  There is a Wood frequency at $k\approx0.921954$.}
\label{fig:GreenConvergenceAngle}
\end{figure}

The main results in this contribution include {\em (i)}~A proof,
presented in Section~\ref{sec:convergence}, establishing that, as the
truncation radius $a$ tends to $+\infty$, the smoothly truncated Green
function converges faster than any negative power of $a$---at least
for arrangements of the period, frequency and incidence angles that
lie away from certain ``Wood configurations'' (for which the Green
function $G^{qper}$ ceases to exist); as well as {\em (ii)}~A new
accelerated integral-equation solver presented in
Section~\ref{sec:integral} which, relying on the aforementioned
windowed Green function, gives rise to a highly-efficient overall
solution method for the problems at hand. In particular,
Theorem~\ref{thm:super-algebraic} below establishes the
super-algebraically fast convergence of the truncated sum to the
three-dimensional quasi-periodic Green function away from Wood
configurations; a corresponding convergence theorem for 1D-periodic
diffraction gratings in $\RR^2$ was presented in~\cite{Monro2007},
cf. also~\cite{BrunoDelourme}. Figures~\ref{fig:GreenConvergenceNormal}
and~\ref{fig:GreenConvergenceAngle} in Section~\ref{sec:convergence}
demonstrate the convergence of the windowed series both near and away
from Wood configurations. The numerical methods presented in
Sections~\ref{sec:integral}, in turn, integrate the windowed Green
function in the context of fast integral-equation
solvers~\cite{br-k1,BrunoKunyansky}. Interestingly, the structure of
the acceleration methodology inherent in these solvers is exploited to
completely avoid evaluation of the windowed Green function at pairs of
surface points, using instead a much smaller number of values of the
Green function on a certain three-dimensional Cartesian grid. A
variety of numerical results presented in Section~\ref{sec:numerical}
demonstrate the character of the resulting solvers for bi-periodic
scattering problems. Green function methods that are valid even at and
around Wood configurations are presented in~\cite{BrunoDelourme} for
two-dimensional configurations, and in Part~II for the
three-dimensional case.

As is well known, for certain wave numbers $k$ and certain Bloch
wave numbers $(\alpha,\beta)$, the lattice sum~\eqref{latticesum1} does
not converge.  This can be seen in the spectral representation of the
Green function that results by applying the Poisson Summation Formula
to the series \eqref{latticesum1}.  Let
$d=\|\vv_1\times\vv_2\|$.  Then
\begin{equation}\label{Fouriersum}
  G^{qper}(\tilde\bfx,z) = \frac{i}{2d} \sum_{j,\ell\in\ZZ}
  \frac{1}{\gamma_{j\ell}} e^{i[(2\pi j\,\vv_1^* \,+\, 2\pi\ell\,\vv_2^*)+\mathbf{k}]\cdot\tilde\bfx}\, e^{i\gamma_{j\ell}|z|}\,,
\end{equation}
in which the propagation constants $\gamma_{j\ell}$ are defined by
\begin{equation}\label{eq:alpha_r}
\vv_{j\ell}^* = (2\pi j\,\vv_1^* \,+\, 2\pi\ell\,\vv_2^*)+\mathbf{k}\,,
\quad
\gamma_{j\ell}=(k^2-\| \vv_{j\ell}^* \|^2)^\frac{1}{2}\,.
\end{equation}
(The branch of the square root that defines $\gamma_{j\ell}$ is
selected in such a way that $\sqrt{1} =1$, and that the branch cut
coincides with the negative imaginary semiaxis.)  The lattice
sum~\eqref{latticesum1} converges if and only if $\gamma_{j\ell}\ne 0$
for all integer pairs $(j,\ell)$.  Configurations for which
$\gamma_{j\ell}$ vanishes for one or more integer pairs $(j,\ell)$ are
known as Wood configurations, or Wood anomalies.
Clearly the expression~(\ref{Fouriersum}) is not meaningful if $\gamma_{j\ell} = 0$ for some integer
pair $(j,\ell)$.  Wood anomalies were first noticed by
Wood~\cite{Wood} and first treated mathematically by
Rayleigh~\cite{Rayleigh}; a brief discussion concerning historical
aspects in these regards can be found in~\cite[Remark
2.2]{BrunoDelourme}.  As shown in~\cite{BrunoDelourme} and Part~II,
Green function methods can still be used at Wood anomalies provided
appropriately defined Green functions are used.

In view of the branch used in equation~\eqref{eq:alpha_r} for the
square root function, Rayleigh waves either decay as $|z|$ increases
(evanescent modes) or are outgoing traveling waves (propagating
modes). Clearly, there exist finitely many propagating modes for any
given configuration.  Wood frequencies are also called ``cutoff
frequencies", since the corresponding Rayleigh wave $e^{i\vv_{j\ell}^*\cdot\tilde\bfx} e^{i\gamma_{j\ell}|z|}$ switches from propagating
to evanescent as the frequency descends below a Wood
value. Rayleigh waves for which $\gamma_{j\ell}$ is small impinge
on the periodic structure at ``grazing incidence'', and they dominate
the sum \eqref{Fouriersum}.  In the limit of a particular combination
of $k$ and $(\alpha, \beta)$, at which one or more $\gamma_{j\ell}$
are zero, the product of the sum multiplied by any one of the
vanishing $\gamma_{j\ell}$'s tends to a $z$-independent linear
combination of exactly grazing waves of the form $e^{i\vv_{j\ell}^*\cdot\tilde\bfx}$.

Challenges in the calculation of the Green function arise from two
main sources, namely
\begin{enumerate}
\item The lattice sum~\eqref{latticesum1} does not converge
  absolutely.  This sum does converge conditionally away from Wood
  anomalies~\cite{BrunoReitich1992}, but its convergence, which
  results from cancellations amongst slowly decreasing terms, is too
  slow to be useful from a computational standpoint.
\item At Wood configurations the lattice sum \eqref{latticesum1} does
  not converge and a denominator in the Rayleigh-wave expansion
  \eqref{Fouriersum} exactly vanishes.  Additionally, the convergence
  of the series~\eqref{latticesum1} increasingly deteriorates as the
  parameters in the problem are varied in such a way that a Wood
  configuration is approached.
\end{enumerate}
The first of these challenges is addressed in the present contribution, and the second is treated  in~\cite{BrunoDelourme}
for the two-dimensional case, and, for three dimensions, in Part
II~\cite{BSTV2}.  

As mentioned above, the proposed approach for summation of the series
is based on smooth windowing of the series~\eqref{latticesum1}.  A similar windowed-summation
technique can be applied to the spectral series~\eqref{Fouriersum}
with similar super-algebraic convergence. A study of the potential
advantages offered by such a strategy is left for future work.

Previous accelerated procedures based on either or both of the spatial
and spectral representations for the Green function $G^{qper}$ give
rise to significantly faster algorithms than does direct summation of
either the expressions~\eqref{latticesum1} or~\eqref{Fouriersum}. The
two-dimensional algorithms (see e.g.~\cite[Section
3.8.2]{Petit-Maystre}) and~\cite{Veysoglu91}) can be perfectly
adequate, but in the three-dimensional context algorithms for
evaluation of quasi-periodic Green functions have remained
inefficient. As a significant reference in these regards we mention
one of the most advanced hybrid approaches previously put forth for
evaluation of periodic Green's functions~\cite{Guerin}, which is based
on use of a combination of spatial and spectral representations as
well as Kummer and Shanks transforms. The hybrid
algorithm~\cite{Guerin} has been reported~\cite{Blezynski}
(cf. also~\cite{Guerin}) to require several milliseconds per
evaluation point. Thus, even for a small discretization consisting of
$N = 6\times 16\times 16$ points (assuming a total of $6$ patches are
used to represent a given scattering surface $S$, and $6\times 6$
discretization points are used in each patch) the number $2 \times
10^{6}$ of evaluations of periodic Green functions which are necessary
to evaluate one matrix-vector product requires a computational time of
at least $2 \times 10^{3}$ seconds. In contrast, as it can be seen in
Table~\ref{errors0}, in the case of periodic two-dimensional arrays of
spheres discretized by means of such a $6\times 16\times 16$ mesh, our
solvers require less than 10 seconds per matrix-vector product (an
improvement factor of a least one-hundred)---and can produce full
scattering results with an error of the order of $10^{-4}$ in a total
of 55 seconds.

As mentioned above, boundary-integral equations based on the proposed
Green-function methods are presented in Section~\ref{sec:integral}.
In particular, Section~\ref{sec:integral} describes the numerical
methods used to implement the proposed fast lattice sums and forward
maps (matrix-vector products) which, upon use of an iterative linear
algebra solver (GMRES) produces the densities in certain
boundary-integral representations of the scattered field. For
definiteness, in all numerical examples it was assumed the scatterers
satisfy sound-soft (Dirichlet) boundary conditions.
Section~\ref{sec:numerical} demonstrates the resulting method by means
of a variety of numerical results. A few concluding remarks are
presented in Section~\ref{sec:concl}.


\section{Proof of fast convergence of smoothly truncated lattice sums}\label{sec:convergence}

Our smooth truncation method proceeds by multiplying the $(m,n)$-th
term of the series \eqref{latticesum1} by the scaled cut-off function
$\chi(\tilde r_{mn}/a)$ defined in equation~\eqref{chi_def}; the smoothly
truncated series is thus given by the finite sum
\begin{equation}\label{windowing}
  G^a(x,y,z) := \frac{1}{4\pi}\sum_{m,n\in\ZZ} \frac{e^{ikr_{mn}}}{r_{mn}}
  e^{-i \mathbf{k}\cdot(m\mathbf{v}_1+n\mathbf{v}_2)} \,\chi\left(\tilde r_{mn}/a\right)\approx\,G^{qper}(x,y,z),
\end{equation}
where of $r_{mn}$ and $\tilde r_{mn}$ are given by (\ref{r})
and~(\ref{rtilde}). The following theorem establishes
the super-algebraic convergence of the truncated lattice sum to the
quasi-periodic Green function for configurations away from Wood
anomalies.

\begin{Theorem}[Windowed Green function at non-Wood frequencies:
  Super-algebraic convergence]\label{thm:super-algebraic}
  Let $\chi(r)$ be an infinitely smooth truncation function which
  equals to $1$ for $r\leq r_1$ and equals $0$ for $r\geq r_2$
  ($0<r_1<r_2$).  If $\gamma_{j\ell}\not=0$ for all
  $(j,\ell)\in\ZZ^2$, then the functions
\begin{equation*}
  G^a(x,y,z) = \frac{1}{4\pi}\sum_{m,n\in\ZZ} \frac{e^{ikr_{mn}}}{r_{mn}}
  e^{-i\mathbf{k}\cdot(m \mathbf{v}_1+n \mathbf{v}_2)} \,\chi\left(\tilde r_{mn}/a\right)
\end{equation*}
converge to the radiating quasi-periodic Green function
$G^\text{qper}(x,y,z)$ super-algebraically fast as the truncation
radius $a$ tends to infinity.  In detail, for each posistive integer
$n$ there exist constants $C_n=C_n(k, \alpha, \beta)$ such that
\begin{equation}\label{superalgebraic}
  | G^a_k(x,y,z) - G^\text{qper}(x,y,z) | < \frac{C_n(k, \alpha, \beta)}{a^n}
\end{equation}
when $a$ is sufficiently large.  The inequality holds uniformly for all
points $(x,y,z)$, excluding the singularities of the Green function
for which $r_{mn}=0$ for some $(m,n)\in\ZZ^2$.  At these points, a
term that is common to $G^a_k$ and $G^\text{qper}$ is infinite.  If
$G^a_k$ and $G^\text{qper}$ are modified by excluding this term then
the correspondingly modified version of
equation~(\ref{superalgebraic}) remains valid.

An analogous estimate holds for $ \| \nabla G^a_k(x,y,z) - \nabla G^\text{qper}(x,y,z) \|$.
\end{Theorem}

\noindent

{\em\bfseries Proof.}  Denote by $\Lambda=\{m\vv_1+n\vv_2 :
m,n\in\ZZ\}$ the lattice of singularities of the Green function,
and denote by $\Lambda^*=\{j\vv_1^*+\ell\vv_2^* : j,\ell\in\ZZ\}$ the dual
lattice.  The dual vectors $\vv_1^*$ and $\vv_2^*$ are defined by $\vv_i^*\cdot\vv_j=\delta_{ij}$.  Initially, we assume that the shift $\tilde\bfx=(x,y)$ from these
positions as well as the Bloch wavenumbers $\alpha$ and $\beta$ are
equal to zero.  Setting $\veps=a^{-1}$, we have for the full and the
truncated sums,
\begin{equation}
 4\pi G^{qper} \,= \sum_{\rr\in\Lambda} \frac{\exp(ik\sqrt{|\rr|^2+z^2\,})}{\sqrt{|\rr|^2+z^2\,}}
\end{equation}
and
\begin{equation}
4\pi G^a \,=  \sum_{\rr\in\Lambda} \chi(\veps |\rr|)  \frac{\exp(ik\sqrt{|\rr|^2+z^2\,})}{\sqrt{|\rr|^2+z^2\,}}\,.
\end{equation}
With the view of utilizing the Poisson summation formula to calculate
the truncated sum we introduce a smooth function $\phi(|\rr|)$ that
vanishes in a neighborhood of $|\rr|=0$ and is equal to $1$ for
$|\rr|\geq r_1$.  For $\veps<1$, the sum is broken into two pieces,
\begin{multline}\label{splitting1}
   \sum_{\rr\in\Lambda} \chi(\veps |\rr|)  \frac{\exp(ik\sqrt{|\rr|^2+z^2\,})}{\sqrt{|\rr|^2+z^2\,}} \\
  =  \sum_{0\not=\rr\in\Lambda} (1-\phi(|\rr|)) \frac{\exp(ik\sqrt{|\rr|^2+z^2\,})}{\sqrt{|\rr|^2+z^2\,}}
  +   \sum_{\rr\in\Lambda} \phi(|\rr|)\chi(\veps|\rr|) \frac{\exp(ik\sqrt{|\rr|^2+z^2\,})}{\sqrt{|\rr|^2+z^2\,}}\,.
\end{multline}
In the first sum on the right,  $\chi$ is omitted as a factor since it equals unity when $\phi\ne 1$. The term is thus independent of the truncation variable $\veps$.
It is easy to check that the fraction in the second term can be
expressed as a product of an exponential function and a Laurent
expansion:
\begin{equation}
  \frac{\exp(ik\sqrt{r^2+z^2})}{\sqrt{r^2+z^2}}
  = \frac{e^{ikr}}{r}g(r)\, , \ \ \ \ \ \ \ \ \  g(r) = 1 + \sum_{j=1}^\infty a_j r^{-j}\,.
\end{equation}
The coefficients $a_j$   are functions of $z$ and the expansion is convergent when $r>|z|$.

We re-express the second sum in (\ref{splitting1}) by means of the
Poisson summation formula:
\begin{equation}\label{PSF1}
  \sum_{\rr\in\Lambda} \phi(|\rr|)\chi(\veps|\rr|)g(|\rr|) \frac{e^{ik|\rr|}}{|\rr|}
  \;=\; \frac{1}{\,d\,} \sum_{\xxi\in\Lambda^*} {\mathcal F} \left[ \phi(|\rr|)\chi(\veps|\rr|)g(|\rr|) \frac{e^{ik|\rr|}}{|\rr|} \right](\xxi),
\end{equation}
where $d=\|\vv_1\times\vv_2\|$.
In what follows we re-express the Fourier transform on the right-hand
side of this equation (which, for brevity, we denote by ${\mathcal
  F}(\xxi)$) in terms of suitable contour integrals.  To do this, we
represent the spatial and Fourier variables in polar coordinates,
$\rr=(r,\theta)$ and $\xxi=(\xi,\gamma)$, and we let $f(r) =
\phi(r)g(r)$, and we thus obtain
\begin{multline}\label{double_integrals}
  {\mathcal F}(\xxi) = \int_0^\infty\!\!\!\int_{-\pi}^\pi f(r)\chi(\veps r) e^{i(k-2\pi\xi\cos(\theta-\gamma))r} d\theta dr  \\
    =\, 2\!\int_0^\infty\! f(r)\chi(\veps r) \!\int_0^\pi e^{i(k-2\pi\xi\cos\theta)r} d\theta dr
    \,=\, 2\!\int_0^\infty\! f(r)\chi(\veps r) \!\int_{-1}^1 e^{i(k-2\pi\xi s)r} \frac{ds}{\sqrt{1-s^2}\,} dr \\
    \,=\, 2\,\int_0^\infty\! f(r)\chi(\veps r)
                       \left(\int_{-1}^{-1-i\infty}e^{i(k-2\pi\xi s)r} \frac{ds}{\sqrt{1-s^2}\,} dr-\int_1^{1-i\infty} e^{i(k-2\pi\xi s)r} \frac{ds}{\sqrt{1-s^2}\,} dr\right).
\end{multline}
The last equality is valid by contour integration in the complex
$s$-plane in view of the exponential decay of the integrand as
$\mathrm{Im}(s)\to-\infty$. We have thus obtained
\begin{equation}\label{FandI}
  {\mathcal F}(\xxi) = 
       2\left(\int_{-1}^{-1-i\infty}I(s)\,  \frac{ds}{\sqrt{1-s^2}\,}-\int_1^{1-i\infty} I(s)\,  \frac{ds}{\sqrt{1-s^2}\,}\right)
\end{equation}
where

\begin{equation}\label{I-expr}
  I(s) = \int_0^\infty f(r)\chi(\veps r) e^{i(k-2\pi\xi s)r} dr\,. 
 \end{equation}
 The integrand~\eqref{I-expr} decays exponentially fast at infinity
 since $\mathrm{Im}(s)<0$. Thus, integration by parts (in which the
 boundary terms vanish because $f(r)$ vanishes near $r=0$) yields
\begin{multline}
 I(s)   = \frac{1}{i(k-2\pi\xi s)} \int_0^\infty \left[ f'(r) + f'(r)(\chi(\veps r)-1) + \veps f(r)\chi'(\veps r) \right]
        e^{i(k-2\pi\xi s)r}dr \\
     = I_0(s) + I_\veps(s)
\end{multline}
where
\begin{equation}
  I_0(s) = \frac{1}{i(k-2\pi\xi s)} \int_0^\infty f'(r) e^{i(k-2\pi\xi s)r}dr
            = \frac{1}{[i(k-2\pi\xi s)]^n} \int_0^\infty f^{(n)}(r) e^{i(k-2\pi\xi s)r}dr,
\end{equation}
and where, noting that $\{\chi = 1\}\supseteq \{\phi\ne 1\}$ for
$\veps<1$, we have $\chi-1=0$, $\chi '=0$ and $f(r)=g(r)$ in the
region $\{\chi = 1\}$, and, thus
\begin{equation}
  I_\veps(s) = \frac{1}{i(k-2\pi\xi s)} \int_0^\infty \left[ g'(r)(\chi(\veps r)-1) + \veps g(r)\chi'(\veps r) \right] e^{i(k-2\pi\xi s)r}dr.
\end{equation}
Thus, introducing a rescaled version $g_\veps$ of the function $g$,
\begin{equation}\label{gepsilon}
  g_\veps(\rho) = g(\rho/\veps) = \frac{\rho}{\sqrt{\rho^2+(\veps z)^2\,}}
  \exp\!\left( \frac{ikz^2\veps}{\rho + \sqrt{\rho^2+(\veps z)^2}\,} \right) 
  = 1 + \sum_{j=1}^\infty a_j \frac{\veps^j}{\rho^j}
\end{equation}
the integrals $I_\varepsilon(s)$  become
\begin{eqnarray*}   
 I_\veps(s) &=& \frac{\veps}{i(k-2\pi\xi s)} \int_0^\infty \left[ g_\veps'(\veps r)(\chi(\veps r)-1) + g_\veps(\veps r)\chi'(\veps r) \right] e^{i(k-2\pi\xi s)r}dr \\
   &=& \frac{1}{i(k-2\pi\xi s)} \int_0^\infty \left[ g_\veps'(\rho)(\chi(\rho)-1) + g_\veps(\rho)\chi'(\rho) \right] e^{i(k-2\pi\xi s)\rho/\veps} d\rho \\
   &=& \frac{(-1)^n\veps^n}{[i(k-2\pi\xi s)]^{n+1}} \int_0^\infty
       \frac{d^n}{d\rho^n} \left[ g_\veps'(\rho)(\chi(\rho)-1) + g_\veps(\rho)\chi'(\rho) \right] e^{i(k-2\pi\xi s)\rho/\veps} d\rho\,.
\end{eqnarray*}
In view of (\ref{FandI}), the splitting $I(s)=I_0(s)+I_\veps(s)$
effects the splitting
\begin{equation}\label{FT_split}
  {\mathcal F}(\xxi) = {\mathcal F}_0(\xxi) + {\mathcal F}_\veps(\xxi),
\end{equation}
for ${\mathcal F}(\xxi)$, where letting
\begin{equation}
  S_\pm(\rho) = 
  \int_{\pm 1}^{\pm 1-i\infty}
  \frac{e^{i(k-2\pi\xi s)\rho/\veps}}{[i(k-2\pi\xi s)]^{n+1}}
  \frac{ds}{\sqrt{1-s^2}}
  =
  -2i \int_0^\infty 
  \frac{e^{(i(k\mp 2\pi\xi)-2\pi\xi t^2)\rho/\veps}}{[i(k\mp2\pi\xi)-2\pi\xi t^2]^{n+1}}
  \frac{dt}{\sqrt{\pm 2i+t^2}\,}
\end{equation}
(the last expression of which incorporates the changes of variables
$s=\pm 1-it^2$) we have denoted
\begin{equation}
  {\mathcal F}_0(\xxi) \,=\, 2\int_0^\infty f^{(n)}(r) 
                                         \left(\int_{-1}^{-1-i\infty}\frac{e^{i(k-2\pi\xi s)r}}{[i(k-2\pi\xi s)]^n}
                                         \frac{ds}{\sqrt{1-s^2}} \,dr\,\!\!-\!\int_1^{1-i\infty}
                                         \frac{e^{i(k-2\pi\xi s)r}}{[i(k-2\pi\xi s)]^n}
                                         \frac{ds}{\sqrt{1-s^2}} \,dr\right)
\end{equation}
and
\begin{equation}\label{twenty_four}
  {\mathcal F}_\veps(\xxi) =
  2\,(-1)^n\veps^n\! \int_0^\infty \frac{d^n}{d\rho^n} \left[
    g_\veps'(\rho)(\chi(\rho)-1) + g_\veps(\rho)\chi'(\rho) \right]
 \left(S_{-}(\rho)-S_{+}(\rho)\right) \,d\rho\,.
\end{equation}

Assume now that $\xi\not=0$; the case $\xi=0$ will be treated
separately.  In view of the hypothesis $k-2\pi\xi\not=0$ the integral
$|S_+|$ admits the finite upper bound
\begin{eqnarray}\label{s_bound_1}
  |S_+(\rho)| &\leq& \int_0^\infty
                       \frac{2\, e^{-2\pi\xi t^2\rho/\veps} dt}{[(k-2\pi\xi)^2+(2\pi\xi t^2)^2]^{\frac{n+1}{2}} (4+t^4)^{\frac{1}{4}}} 
              \,\leq\,  \int_0^\infty \frac{\sqrt{2}\, e^{-2\pi\xi t^2\rho/\veps} dt}{|k-2\pi\xi|^{n+1}} \\
     &=& \left( \frac{\veps}{\rho\,\xi} \right)^{\frac{1}{2}} \frac{1}{|k-2\pi\xi|^{n+1}}
             \int_0^\infty e^{-\pi t^2} dt
                \,=\,  \frac{1}{2} \left( \frac{\veps}{\rho\,\xi} \right)^{\frac{1}{2}}
                     \frac{1}{|k-2\pi\xi|^{n+1}} \,.
\end{eqnarray}
Analogously, in view of the assumption $k+2\pi\xi\not=0$ we obtain
\begin{equation}\label{s_bound_2}
  |S_-(\rho)| \,\leq\, \frac{1}{2} \left( \frac{\veps}{\rho\,\xi} \right)^{\frac{1}{2}}
                     \frac{1}{|k+2\pi\xi|^{n+1}} \,.
\end{equation}
Returning to the expression for ${\mathcal F}_\veps(\xxi)$ above,
observe that, since $\chi(\rho)=1$ for $\rho\leq r_1$, and
$\chi(\rho)=0$ for $\rho\geq r_2$, the integral in $\rho$ from $0$ to
$\infty$ in~\eqref{twenty_four} can be re-expressed in the form
${\mathcal F}_\veps(\xxi) = {\mathcal F}^1_\veps(\xxi) + {\mathcal
  F}^2_\veps(\xxi)$ where
\begin{eqnarray*}
  {\mathcal F}^1_\veps(\xxi) &=&
  2(-1)^n\veps^n\!\! \int_{r_1}^{r_2}
               \frac{d^n}{d\rho^n} \left[ g_\veps'(\rho)(\chi(\rho)-1) + g_\veps(\rho)\chi'(\rho) \right]
                                         \big(S_-(\rho) - S_+(\rho)\big) \,d\rho\,\quad\mbox{and} \\
   {\mathcal F}^2_\veps(\xxi) &=&
   2(-1)^{n+1}\veps^n \int_{r_2}^\infty g_\veps^{(n+1)}(\rho) \big( S_-(\rho) - S_+(\rho) \big) d\rho\,.
\end{eqnarray*}
The bounds~\eqref{s_bound_1} and~\eqref{s_bound_2} thus imply
\begin{equation}\label{f_eps_bd_1}
  \left| {\mathcal F}^1_\veps(\xxi) \right| \leq
     \frac{\veps^{n+\frac{1}{2}}}{\xi^\frac{1}{2}}
     \left( \frac{1}{\left| k-2\pi\xi \right|^{n+1}} + \frac{1}{\left| k+2\pi\xi \right|^{n+1}} \right)
     \int_{r_1}^{r_2} \left| \frac{d^n}{d\rho^n} \left[ g_\veps'(\rho)(\chi(\rho)-1) + g_\veps(\rho)\chi'(\rho) \right] \right| \frac{d\rho}{\rho^{\frac{1}{2}}}\,.
\end{equation}

Clearly, as $\veps\to0$ the functions $g_\veps(\rho)$ converge to $1$
uniformly over the interval $[r_1,r_2]$, and thus the
integral~\eqref{f_eps_bd_1} integral converges to $\int_{r_1}^{r_2}
\chi^{(n+1)}(\rho) \rho^{-1/2} d\rho$ in this limit. In particular
these integrals are bounded by a constant $C^1_n>0$ for all $\veps<1$
and we have
\begin{equation}\label{Feps1bound}
  \left| {\mathcal F}^1_\veps(\xxi) \right| \leq
  C^1_n \frac{\veps^{n+\frac{1}{2}}}{\xi^\frac{1}{2}}
     \left( \frac{1}{\left| k-2\pi\xi \right|^{n+1}} + \frac{1}{\left| k+2\pi\xi \right|^{n+1}} \right).
\end{equation}
Similarly, for ${\mathcal F}^2_\veps(\xxi)$ we have
\begin{equation}
  \left| {\mathcal F}^2_\veps(\xxi) \right| \leq
     \frac{\veps^{n+\frac{1}{2}}}{\xi^\frac{1}{2}} 
     \left( \frac{1}{\left| k-2\pi\xi \right|^{n+1}} + \frac{1}{\left| k+2\pi\xi \right|^{n+1}} \right)
     \int_{r_2}^\infty \left| g_\veps^{(n+1)}(\rho) \right| \frac{d\rho}{\rho^\frac{1}{2}}\,.
\end{equation}
But from~\eqref{gepsilon} we obtain
\begin{equation}
  g_\veps^{(n+1)}(\rho) = \frac{(-1)^{n+1}}{\rho^{n+1}}
                     \sum_{j=1}^\infty a_j \frac{(j+n)!}{(j_1)!} \frac{\veps^j}{\rho^j}\, ,
\end{equation}
and, we thus see that, for $\veps$ sufficiently small,
$\int_{r_2}^\infty \big| g_\veps^{(n+1)}(\rho)
\big|{\rho^{-\frac{1}{2}}} d\rho$ is bounded by a certain constant
$C^2_n$, so that
\begin{equation}\label{Feps2bound}
  \left| {\mathcal F}^2_\veps(\xxi) \right| \leq
  C^2_n \frac{\veps^{n+\frac{1}{2}}}{\xi^\frac{1}{2}}
     \left( \frac{1}{\left| k-2\pi\xi \right|^{n+1}} + \frac{1}{\left| k+2\pi\xi \right|^{n+1}} \right).
\end{equation}
Combining the estimates ${\mathcal F}^1_\veps(\xxi)$ and ${\mathcal
  F}^2_\veps(\xxi)$ we thus find that there exists a constant $C^3_n$
such that
\begin{equation}\label{Fepsbound}
  \left| {\mathcal F}_\veps(\xxi) \right| \leq
  \veps^{n+\frac{1}{2}} \, \frac{C^3_n}{\xi^\frac{1}{2}}
     \left( \frac{1}{\left| k-2\pi\xi \right|^{n+1}} + \frac{1}{\left| k+2\pi\xi \right|^{n+1}} \right).
\end{equation}

For $\xi=0$, in turn, we have
\begin{multline}\label{xi_equal_zero}
  {\mathcal F}(\mathbf{0}) = \int_0^\infty\int_{-\pi}^\pi f(r)\chi(\veps r) e^{ikr} d\theta dr
       = 2\pi \int_0^\infty f(r)\chi(\veps r) e^{ikr} dr = \\
   = -\frac{2\pi}{ik} \int_0^\infty f'(r) e^{ikr} dr
           -\frac{2\pi}{ik} \int_0^\infty \left[ g_\veps'(\rho)(\chi(\rho)-1) + g_\veps(\rho) \chi'(\rho) \right]e^{ik\rho/\veps} d\rho \\
   = -\frac{2\pi}{ik} \int_0^\infty f'(r) e^{ikr} dr
   + (-1)^n \frac{2\pi}{ik} \left( \frac{\veps}{ik} \right)^{n+1}
       \int_0^\infty \frac{d^n}{d\rho^n} \left[ g_\veps'(\rho)(\chi(\rho)-1) + g_\veps(\rho) \chi'(\rho) \right]e^{ik\rho/\veps} d\rho \\
    = {\mathcal F}_0(\mathbf{0}) + {\mathcal F}_\veps(\mathbf{0})\,.
\end{multline}
Again, ${\mathcal F}_0(\mathbf{0})$ is independent of $\veps$ and the integral in ${\mathcal F}_\veps(\mathbf{0})$ has a limit as $\veps\to0$.  Thus one obtains constants $C_n^0>0$ such that $|{\mathcal F}_\veps(\mathbf{0})|\leq C_n^0 \veps^{n+1}$.

The estimates above now allow us to now establish the convergence as
$\veps\to 0$ of the series on the right-hand side of
equation~(\ref{PSF1}).  If $n\geq1$, then as long as
$2\pi|\xxi|\not=|k|$ for all $\xxi\in\Lambda^*$, the sum of $\left|
  {\mathcal F}_\veps(\xxi) \right|$ over all $\xxi\in\Lambda^*$ is
convergent, and one obtains
\begin{equation}\label{Fepsilonbound}
  \sum_{\xxi\in\Lambda^*} \left|{\mathcal F}_\veps(\xxi) \right| \leq C_n\, \veps^{n+\frac{1}{2}}\,.
\end{equation}
The Poisson Summation Formula now gives
\begin{equation}\label{PSF2}
  d\, \sum_{\rr\in\Lambda} \phi(|\rr|)\chi(\veps|\rr|)g(|\rr|) \frac{e^{ik|\rr|}}{|\rr|}
  = \sum_{\xxi\in\Lambda^*} {\mathcal F}_0 (\xxi) + \sum_{\xxi\in\Lambda^*} {\mathcal F}_\veps (\xxi).
\end{equation}
But the first term on the right hand side of this equation is
independent of $\veps$, and, in view of~(\ref{Fepsilonbound}), the
second term on the right hand side tends to zero super-algebraically
fast. It follows that the sum on the left hand side of~\eqref{PSF2}
converges super-algebraically fast, as needed.

\vspace{0.1cm}
\noindent
{\bf Inclusion of the Bloch quasi-periodicity factors} in the lattice sum
can now be accomplished by replacing the expression
$\exp(ik\sqrt{|\rr|^2+z^2})/\sqrt{|\rr|^2+z^2}$ by
\begin{equation}
  \frac{\exp(ik\sqrt{|\rr|^2+z^2\,})}{\sqrt{|\rr|^2+z^2\,}} e^{-i{\mathbf{k}}\cdot\rr},
\end{equation}
where ${\mathbf{k}}=\alpha\vv_1^*+\beta\vv_2^*$.
Equation (\ref{PSF2}) becomes
\begin{equation}\label{PSF3}
  d\, \sum_{\rr\in\Lambda} \phi(|\rr|)\chi(\veps|\rr|)g(|\rr|) \frac{e^{ik|\rr|}}{|\rr|}e^{-i{\mathbf{k}}\cdot\rr}
  = \sum_{\xxi\in\Lambda^*} {\mathcal F}_0 \left(\xxi+\frac{{\mathbf{k}}}{2\pi}\right) + \sum_{\xxi\in\Lambda^*} {\mathcal F}_\veps \left(\xxi+\frac{{\mathbf{k}}}{2\pi}\right),
\end{equation}
The bound (\ref{Fepsbound}), shifted by ${\mathbf{k}}/2\pi$, is
\begin{equation*}
  \left| {\mathcal F}_\veps \left(\xxi+\frac{{\mathbf{k}}}{2\pi}\right) \right| \leq
  \veps^{n+\frac{1}{2}} \, \frac{C^3_n}{\xi^\frac{1}{2}}
     \left( \frac{1}{\big| k-2\pi\left| \xxi + {\mathbf{k}}/(2\pi) \right| \big|^{n+1}} + \frac{1}{\big| k+2\pi\left| \xxi + {\mathbf{k}}/(2\pi) \right| \big|^{n+1}} \right),
\end{equation*}
which is valid whenever
\begin{equation}\label{Wood2}
  k^2 \not= \left| 2\pi \xxi + {\mathbf{k}} \right|^2\,.
\end{equation}
The validity of (\ref{Wood2}) for all $\xi\in\ZZ^2$ is exactly the condition that $(k,\alpha,\beta)$ is not a Wood triple.
 
\vspace{0.1cm}
\noindent {\bf Inclusion of a shift in $\rr$} by a fixed vector
$\rr'=(x,y)$. Consider the lattice sum of the quantities
\begin{equation}
    \frac{\exp(ik\sqrt{|\rr - \rr'|^2+z^2\,})}{\sqrt{|\rr - \rr'|^2+z^2\,}},
\end{equation}
in which we have taken ${\mathbf{k}}=0$.  The case ${\mathbf{k}}\not=0$ is again treated by shifting the Fourier variable $\xxi$ as shown above.
As the cutoff functions $\phi$ and $\chi$ are also shifted, there ensues a mere exponential factor in the Fourier transform, and equation (\ref{PSF2}) becomes
\begin{equation}\label{PSF4}
  d\, \sum_{\rr\in\Lambda} \phi(|\rr - \rr'|)\chi(\veps|\rr - \rr'|)g(|\rr - \rr'|) \frac{e^{ik|\rr - \rr'|}}{|\rr - \rr'|}
  \,=\, \sum_{\xxi\in\Lambda^*} {\mathcal F}_0 (\xxi) e^{-2\pi i \xxi\cdot\rr'}
  + \sum_{\xxi\in\Lambda^*} {\mathcal F}_\veps (\xxi) e^{-2\pi i \xxi\cdot\rr'}.
\end{equation}
The bound (\ref{Fepsilonbound}) persists,
\begin{equation}
  \bigg| \sum_{\xxi\in\Lambda^*} {\mathcal F}_\veps(\xxi) e^{-2\pi i \xxi\cdot\rr'} \bigg| \leq C_n\, \veps^{n+\frac{1}{2}}\,,
\end{equation}
and one again obtains super-algebraic convergence.

\vspace{0.1cm}
\noindent
{\bf Error bound for the gradient of the Green function.}
The gradient of the monopole $e^{ikr}/r$ is given by the equations 
 \begin{equation}
\label{partial_x_or_y}
\frac{\partial}{\partial x}\frac{e^{ikr}}{r}=\left(ik\cos \theta-\frac{\cos \theta}{r}\right)\frac{e^{ikr}}{r}, \ \ \ \ \ \ \
\frac{\partial}{\partial y}\frac{e^{ikr}}{r}=\left(ik\sin \theta-\frac{\sin \theta}{r}\right)\frac{e^{ikr}}{r}.
\end{equation}
\begin{equation}
\label{partialz}
 \frac{\partial}{\partial z}\frac{e^{ikr}}{r}
=\left(ik\frac{z}{r}-\frac{z}{r^2}\right)\frac{e^{ikr}}{r}. 
\end{equation}
It suffices to show that   the error bound proven in the theorem remains true, if the monopole $\frac{e^{ikr}}{r}$ is replaced in the proof by any of the terms of the above equations. These terms  are products of  the monopole multiplied  by $r^{-1}$, or by $\cos\theta$, or by $\sin\theta$, or by a selection of two of these factors.  The bound is clearly preserved when  multiplying the monopole by $z$, since the latter factors out of the summation that constitutes the Green function. 

Multiplying the monopole by $r^{-1}$ or by $r^{-2}$  corresponds to introducing the factor $\veps/\rho$ or $(\veps/ \rho)^2$ respectively in the subsequent  integration over $\rho$, 
 thus enhancing the error bound by one or two orders  in $\veps$. The integrand is zero  (see explanation following \eqref{twenty_four}) when $\rho<r_1$, thus the denominator $\rho$ is no  cause of concern.

The following observations show that the error bound is preserved in the terms in which the monopole is multiplied by \,$\cos\theta$\, or by \,$\sin\theta$. 
\begin{itemize}
 \item The  first double integral of \eqref{double_integrals} acquires the factors $\cos\theta$ or  $\sin\theta$ in its integrand.   
Thus, the second double integral in \eqref{double_integrals} (obtained by the 
change of the integration variable $\theta\to\theta+\gamma$) exhibits the factors   $\cos(\theta+\gamma)$ or $\sin(\theta+\gamma)$ that can be split into a linear combination of  $\cos\theta$ and $\sin\theta$, with the corresponding splitting of the integral. 
\item The double integral that contains  the factor  $\sin\theta$ is  equal to zero; the integrand of the integration with respect to $\theta$ is an exact derivative and the integration is over the closed loop from $-\pi$ to $\pi$.  
 \item What is left is the second double integral in \eqref{double_integrals} with the extra factor  $\cos\theta$ in the integrand. The change of the variable of integration  $\cos\theta=s$ leads to having an extra factor $s$ in the subsequent integrals with respect to $s$.  This introduces the extra factor $|s|=|\pm 1-it^2|$ into the numerator of the first integral in \eqref{s_bound_1}. 
\item Following the  change of variable $s=\pm1-it^2$, the factor $|s|$ is replaced by its upper bound  $1+t^2$ 
and the integral is split accordingly into a sum of two integrals. The first integral is exactly the one that provides the error bound of the theorem.  The extra factor $t^2$ in the second  integral provides the extra factor $\frac{\veps}{\xi\rho}$ in the bounds \eqref{s_bound_1} and  \eqref{s_bound_2} when $\xi\ne 0$.
Thus, the error bound of the theorem is preserved in this case. 
\item If $\xi=0$, the first integral in \eqref{xi_equal_zero}  has the factor \,$\cos\theta$\, or \,$\sin\theta$\, that integrates to zero. 
\end{itemize}

\hfill$\blacksquare$


\section {Fast high-order integral solvers for problems of scattering by doubly-periodic structures}\label{sec:integral}

For definitness we restrict our treatment to diffractive structures
consisting of arrays of separated obstacles arranged in a
two-dimensional periodic fashion in three-dimensional space.  Thus,
denoting by $\Omega\subset \mathbb{R}^3$ an open connected set the
region occupied by a ``reference obstacle'' (which could itself be
given by the union of a number of connected components) and letting
$S=\partial \Omega$ denote its boundary (the reference scattering
boundary), the overall three-dimensional bi-periodic scattering
structure and its boundary are given by 
\begin{equation}\label{omega_per}
  \Omega_{per}=\bigcup_{(m,n)\in\mathbb{Z}\times\mathbb{Z}}\Omega_{m,n}\;,\quad S_{per}=\bigcup_{(m,n)\in\mathbb{Z}\times\mathbb{Z}}S_{m,n}\;,
\end{equation}
respectively, where, we have set $\Omega_{m,n}=\Omega
-m\mathbf{v}_1-n\mathbf{v}_2$ and $S_{m,n}=S
-m\mathbf{v}_1-n\mathbf{v}_2$, $m,n\in\mathbb{Z}$.
It will be assumed now that the sets $\Omega_{m,n}$, as well as
their boundaries, are pairwise disjoint. Consider the 
sound-soft scattering problem 
\begin{eqnarray}\label{eq:Helmholtz1}
\Delta u+k^2 u&=&0\ {\rm in}\ \mathbb{R}^3\setminus \Omega_{per}\nonumber\\
u&=&-u^{inc}\ {\rm on}\ \partial\Omega_{per}, 
\end{eqnarray} 
in which an incident plane wave 
\begin{equation}\label{eq:planewave}
  u^{inc}(\mathbf x)= \exp(ik\mathbf{d}\cdot\mathbf{x})=\exp[i(\mathbf{k}\cdot\tilde\bfx-\gamma z)],
\end{equation}
with $|\mathbf{k}|^2+\gamma^2=k^2$,
illuminates the structure from above and thus gives rise to a
scattered field $u$.
Owing to the periodicity of the domain $\Omega_{per}$, in the regions
$\Omega_+$ and $\Omega_-$ above and below the array
($\Omega_+=\{\mathbf{x}: z>\max{z'},(x',y',z')\in\Omega_{per}\}$, and
$\Omega_-=\{\mathbf{x}:z<\min{z'},(x',y',z')\in \Omega_{per}\}$) the
fields satisfy radiation conditions expressed in terms of the
classical Rayleigh expansions: the scattered fields $u^+$ and $u^-$ in
the respective regions $\Omega_+$ and $\Omega_-$ must be ``outgoing'',
that is, they must admit Rayleigh expansions of the form
\begin{eqnarray}
u^+(\mathbf x) &=& \sum_{j,\ell\in\ZZ}B^+_{j\ell}\, \exp[i(2\pi j\mathbf{v}_1^{*}+2\pi \ell\mathbf{v}_2^{*}+\mathbf{k})\cdot\tilde\bfx]\exp[i\gamma_{j\ell}z],\quad \mathbf{x}\in\Omega_+\label{eq:rad_cond1}\\
u^-(\mathbf x) &=&\sum_{j,\ell\in\ZZ} B^-_{j\ell}\, \exp[i(2\pi j\mathbf{v}_1^{*}+2\pi \ell\mathbf{v}_2^{*}+\mathbf{k})\cdot\tilde\bfx]\exp[-i\gamma_{j\ell}z],\quad \mathbf{x}\in\Omega_-.\label{eq:rad_cond2}
\end{eqnarray}
wherein no waves in $\Omega_+$ propagate downwards, and no waves in
$\Omega_-$ propagate upwards.  

Using the outgoing free-space Green function $G_k(\mathbf
z)=\frac{e^{ik|\mathbf z|}}{4\pi|\mathbf z|}$, the scattered field $u$
is sought in the form of a combined-field layer potential
\begin{equation}\label{eq:SLf1}
u(\mathbf x)=\int_{S_{per}}\frac{\partial G_k(\mathbf x-\mathbf x')}{\partial \mathbf{n}(\mathbf{x'})}\varphi_\textit{qper}(\mathbf x')ds(\mathbf x') +i\eta\int_{S_{per}}G_k(\mathbf x-\mathbf x')\varphi_\textit{qper}(\mathbf x')ds(\mathbf x')
\end{equation}
with unknown surface density $\varphi_\textit{qper}$. Here
$\mathbf{n}$ is the outer unit normal to $S_{per}$ and
$\eta\in\mathbb{R}$ denotes a coupling constant.  The unknown density
$\varphi_\textit{qper}$ is the solution of the combined field integral
equation
\begin{eqnarray}\label{eq:SLfE1}
\frac{1}{2}\varphi_\textit{qper}(\mathbf{x})+\int_{S_{per}}\frac{\partial G_k(\mathbf x-\mathbf x')}{\partial \mathbf{n}(\mathbf{x'})}\varphi_\textit{qper}(\mathbf x')ds(\mathbf x') &+&i\eta \int_{S_{per}}G_k(\mathbf x-\mathbf x')\varphi_\textit{qper}(\mathbf x')ds(\mathbf x')\nonumber\\
&=&-\exp(ik\mathbf{d}\cdot\mathbf{x}),\ \mathbf{x}\in S_{per}
\end{eqnarray}
which enforces the sound-soft boundary condition under consideration.
The well known term $\frac{1}{2}\varphi_\textit{qper}$
in~\eqref{eq:SLfE1} arises as a singular contribution of the first
integral in~\eqref{eq:SLf1} in the limit as $\mathbf{x}$ approaches
the boundary.

Equations~(\ref{eq:SLfE1}) can be rewritten in a form that involves integration {\it over the
  reference boundary $S$ only}.  The corresponding integral equations
make use of the $(\alpha,\beta)$-quasi-periodic Green function~\eqref{latticesum1}, in which $\bfx$ is replaced by the difference
between source and influence points,
\begin{equation}\label{eq:Greenp}
G_k^\textit{qper}(\mathbf{x}-\mathbf{x}')=\sum_{m,n=-\infty}^\infty G_k((x-x,y-y')+m\mathbf{v}_1+n\mathbf{v}_2,z-z')e^{-i m\mathbf{k}\cdot\mathbf{v}_1}e^{-i n\mathbf{k}\cdot\mathbf{v}_2}. 
\end{equation}
The integral equation~\eqref{eq:SLfE1} can equivalently be expressed in the
form
\begin{eqnarray}\label{eq:SLfE1_2}
  \frac{1}{2}\varphi_\textit{qper}(\mathbf{x})+\sum_{m,n\in\ZZ}\int_{S_{m,n}}\mathcal{G}(\mathbf x-\mathbf x')\varphi_\textit{qper}(\mathbf x')ds(\mathbf x')= -\exp(ik\mathbf{d}\cdot\mathbf{x}),\ \mathbf{x}\in S_{per},
\end{eqnarray}
where 
\begin{equation}
\mathcal{G}(\mathbf x-\mathbf x') = \frac{\partial G_k(\mathbf
  x-\mathbf x')}{\partial \mathbf{n}(\mathbf{x'})} + i\eta
G_k(\mathbf x-\mathbf x').
\end{equation}

Denoting by $\varphi$ the restriction of $\varphi_\textit{qper}$ to the
reference boundary $S$ and taking into account the quasi-periodicity
of the density $\varphi_\textit{qper}$, the integral equation~\eqref{eq:SLfE1}
can be re-expressed in the form
\begin{eqnarray}\label{eq:SLfEPer0}
\frac{\varphi(\mathbf{x})}{2}+\int_{S}\frac{\partial G_k^{per}(\mathbf x-\mathbf{x}')}{\partial \mathbf{n}(\mathbf{x'})}\varphi(\mathbf{x}')ds(\mathbf{x}')&+& i\eta \int_{S}G_k^{per}(\mathbf x-\mathbf{x}')\varphi(\mathbf{x}')ds(\mathbf{x}')\nonumber\\
&=&-\exp(ik\mathbf{d}\cdot\mathbf{x}),\ \mathbf{x}\in S.
\end{eqnarray}
Thus, solution of either equation~\eqref{eq:SLfE1_2}
or~\eqref{eq:SLfEPer0} produces the density $\varphi(\mathbf x)$
which, upon insertion into~\eqref{eq:SLf1} gives rise to the desired
quasi-periodic scattered field. Note that, in view of its
quasi-periodicity, the unknown $\varphi$ is determined throughout
$S_{per}$ by its values on the unit cell $S$---and, thus, testing on
$S$ should suffice to determine $\varphi$ uniquely. Indeed, the
uniqueness of the problem thus posed, which is not pursued here at any
length, can be established by using the periodic Green function as in
equation~\eqref{eq:SLfEPer0} together with a proof similar to the one
for the bounded obstacle case~\cite{ColtonKress1}.


\subsection{High-order evaluation of quasi-periodic layer potentials}\label{sec:evaluation}

Our Nystr\"om approach relies on use of high-order quadratures for
evaluation of the integral operators
\begin{equation*}
  (\mathcal{K}_{m,n}\varphi)(\mathbf{x})=\int_{S_{m,n}} \mathcal{G}(\mathbf{x}-\mathbf{x}')\varphi(\mathbf{x'})ds(\mathbf{x'})
\end{equation*}
in equation~\eqref{eq:SLfE1_2} for $\mathbf{x}\in S$, where
$\varphi=\varphi_\textit{qper}$ is a quasi-periodic integral density
defined on $S_{per}$; as noted in the previous section, testing (and
thus operation evaluation) for $\mathbf{x}\in S$ suffices to determine
the solution $\varphi$. Once such operators have been discretized and
evaluated numerically for a given quasi-periodic function $\varphi$
the solution of the problem can be obtained by means of an iterative
linear algebra solver such a s GMRES~\cite{SaadSchultz}.

We first consider a quadrature algorithm for the operator
$\mathcal{K}=\mathcal{K}_{0,0}$, which is given by
\begin{equation}
\label{eq::compIO}
(\mathcal{K} \varphi)(\mathbf x) = \int_S {\mathcal G}(\mathbf
x-\mathbf{x}')\varphi(\mathbf{x}')ds(\mathbf{x'}),\quad \mathbf{x}\in S.
\end{equation}
We note that this integral operator coincides with the one introduced
in~\cite{br-k1} for the problem of acoustic scattering by a {\em
  bounded obstacle} $S$ under sound-soft boundary conditions. In fact,
the algorithm we propose for evaluation of the integral operators
in~\eqref{eq:SLfEPer0} results as an outgrowth of the fast high-order
methods presented in that reference. (Extensions of these methods to
sound-hard and electromagnetic problems can be found in~\cite{bept}
and~\cite{bet}.) Thus, in order to convey the main ideas underlying
our periodic-structure solver, we first briefly review the
algorithm~\cite{br-k1}.

The bounded-scatterer algorithm~\cite{br-k1} evaluates the integral
operator $\mathcal{K}$ in two stages, namely (a) Evaluation of the
adjacent/singular interactions (i.e. integration for $\mathbf{x'}$ in
areas close to $\mathbf{x}$), and (b) Accelerated evaluation of
nonadjacent interactions (that is, accelerated integration for
$\mathbf{x'}$ away from $\mathbf{x}$). The decomposition into
adjancent and non-adjacent contributions is effected in this method by
means of {\em floating partitions of unity}---that is, pairs of
functions of the form
$(\eta_\mathbf{x}(\mathbf{x'}),1-\eta_\mathbf{x}(\mathbf {x'}))$,
where $\eta_\mathbf{x}$ is a windowing function with a ``small''
support, which equals 1 in a neighborhood of
$\mathbf{x}$. Additionally, the approach~\cite{br-k1} relies on use of
smooth parametrizations of the surface $S$ via a family of overlapping
two-dimensional parameter patches $\mathcal{P}^{\ell},\ell=1,\ldots P$
along with smooth mappings $\mathcal{P}^{\ell}$ from parameter sets
$\mathcal{H}^\ell$ in two-dimensional space (where actual integrations
are performed), as well as partitions of unity subordinated to the
overlapping patch decomposition of the surface.  i.e smooth functions
$w_\ell$ supported on $\mathcal{P}^\ell$ such that $\sum_\ell
w_\ell=1$ throughout $S$. This framework allows us to reduce the
integration of the density $\varphi$ over the surface $S$ to
integration of smooth functions $\varphi^\ell$ compactly supported in
the planar sets $\mathcal{H}^\ell$. The latter calculations require
analytic resolution of weakly singular Green's functions (i.e. the
order of the singularity is $\mathcal O(|\mathbf x
-\mathbf{x'}|^{-1})$) which is performed via polar changes of
variables (whose Jacobian cancels the Green-function singularity)
together with interpolation procedures that facilitate evaluations of
the surface density at radial integration points~\cite{br-k1}.

\subsection{Reference acceleration cell\label{fr_space_accel}}

As indicated at the beginning of Section~\ref{sec:integral}, in this
paper we consider bi-periodic structures of the
form~\eqref{omega_per}. The present Section~\ref{fr_space_accel}
constructs a certain ``reference acceleration cell'' (associated with
the ``reference domain'' $\Omega = \Omega_{0,0}$) which equals a cubic
domain $C$ of side $A$ that contains $\Omega$. The cell $C$ is
equipped with a certain acceleration infrastructure which is based on
a corresponding acceleration technique introduced in~\cite{br-k1}. In
fact, the reference acceleration cell will be utilized as an element
in a method for FFT acceleration for the problem of scattering by the
complete periodic structure $\Omega_\textit{per}$. Here and through
the end of Section~\ref{sec:integral} the presentation assumes a
degree of familiarity with the acceleration methodology presented in
reference~\cite{br-k1}.

The acceleration infrastructure presented in that reference, which is
designed to enable efficient FFT-based acceleration for the numerical
evaluation of the integral operator
\begin{eqnarray}\label{S_oper}
  \int_{S}\mathcal{G}(\mathbf x-\mathbf x')\varphi(\mathbf x')ds(\mathbf x'),\quad{\mathbf x}\in S, 
\end{eqnarray}
(the term $m=n=0$ in~\eqref{eq:SLfE1_2} restricted to ${\mathbf x}\in
S$) proceeds at first by partitioning the cube $C$ into a number
$L^{3}$ of identical cubic cells $c_i$, where $L$ denotes an integer.
The pairs $(A,L)$ of parameters must be adjusted, if necessary, in
order to ensure that the cells $c_i$ do not admit inner acoustic
resonances (eigenfunctions of the Laplace operator with homogeneous
Dirichlet boundary conditions).

The acceleration algorithm~\cite{br-k1} then constructs approximations
which are obtained by substitution of the surface ``true'' sources
within $c_i$ (or, more precisely, of the fields that result from
discrete integration of the product of the kernel $\mathcal G$ and the
density $\varphi$ for all discretization points within $c_i$) by
``equivalent sources'' on a set $\Pi^\ell_i$ ($\ell=1,2,3$) which
equals the union of a pair of parallel circular domains which contain
the faces of $c_i$ that are parallel to the plane $x_\ell=0$. Clearly,
there are three different such approximations. In all three cases the
acoustic fields generated by the $c_i$-equivalent sources approximate
with high order accuracy the fields produced by the true $c_i$ sources
at all cells $c_j$ non-adjacent to $c_i$. The precise concept of
adjacency in~\cite{br-k1} results from a requirement that the
approximation corresponding to a given cell $c_i$ be valid, with
exponentially small errors, outside a concentric cube $\mathcal{S}_i$
of side three times larger than that of $c_i$.  For efficiency the
method relies on use of equivalent sources (acoustic monopoles and
dipoles) as described in what follows. For a given integral density,
and for each cell $c_i$ a set of equivalent sources (acoustic
monopoles $\xi_{i,j}^{(m)\ell}\ G_k(\mathbf x-\mathbf x_{i,j}^{\ell})$
and dipoles $\xi_{i,j}^{(d)\ell}\ \partial G_k(\mathbf x-\mathbf
x_{i,j}^{\ell})/ \partial x_{\ell}$) placed at points
$\mathbf{x}_{i,j}^{\ell},j=1,\cdots,M^{equiv}$ contained within the
union of two circular domains concentric with and circumscribing the
faces of $c_i$, whose radii are selected in accordance with the
prescriptions in {\cite{br-k1}}.  The fields $\psi^{c_i,true}$
radiated by the $c_i$-true sources are approximated by fields
$\psi^{c_{i},eq}$ radiated by the $c_i$ equivalent sources
\begin{equation}
\label{eq:field_ac_eq_sources}
\psi_{0,0}^{c_i,eq}(\mathbf x)=\sum_{j=1}^{\frac{1}{2}M^{equiv}}\left(\xi_{i,j}^{(m)\ell}\ G_k(\mathbf x-\mathbf x_{i,j}^{\ell})+\xi_{i,j}^{(d)\ell}\frac{\partial G_k(\mathbf x-\mathbf x_{i,j}^{\ell})}{\partial x_{\ell}}\right),\ \mathbf{x}\not\in \mathcal{S}_i.  
\end{equation}
For a given number $M^{equiv}$ of equivalent sources (selected so as
to maintain a given accuracy), the unknown monopole and dipole
intensities in (\ref{eq:field_ac_eq_sources}) are chosen so as to
minimize in the mean-square norm the differences
$(\psi^{c_i,eq}(\mathbf{x})-\psi^{c_i,true}(\mathbf x))$ as $\mathbf
x$ varies over a number $n^{coll}$ collocation points on $\partial
\mathcal{S}_i$. Hence, the intensities in
(\ref{eq:field_ac_eq_sources}) are obtained in practice as the
least-squares solution of an overdetermined linear system
$\mathbf{A}\xi=\mathbf{b}$ where $\mathbf{A}$ is an $n^{coll}\times
M^{equiv}$ matrix. As discussed in Sections~\ref{corr}
and~\ref{adj_int} below the method is completed via a sequence of
steps which include 1)~FFTs (which are used to evaluate the Cartesian
convolutions that result from use of equivalent sources);
2)~Correction of certain errors that arise per step~1), which are
inevitable in the FFT-based operation of convolution with the Green
function, and which result from ``incorrect'' use of equivalent
sources for near interactions; and finally, 3)~High-order evaluation
of surface values from the values at the FFT grid. But, before such
discussions we consider certain specializations of the methods above
to the periodic context which, in conjunction with the windowing
methodology used in this paper, have proven specially efficient.
\subsubsection{Green-function contributions from periodic translates
  of the reference cell}
It is easy to check that the set of equivalent sources for the
reference scatterer $\Omega$, as computed per the methodology
described in Section~\ref{fr_space_accel}, can be utilized to
produce---by means of simple algegraic manipulations---the
corresponding equivalent sources for any periodic translation of the
unit-cell. Indeed, denoting by
$(\mathcal{K}_{mn}\varphi_\textit{qper})(\mathbf{x})$ the $(m,n)$-th
term on the left-hand sum in equation~\eqref{eq:SLfE1_2} and since for
$\mathbf{x}\in S$ we have
$\varphi_\textit{qper}(\mathbf{x}-m\mathbf{v}_1-n\mathbf{v}_2)= e^{-i(m \mathbf{k}\cdot\mathbf{v}_1+
  n \mathbf{k}\cdot\mathbf{v}_2)}\varphi (\mathbf{x})$, it follows that, for $\mathbf{x}\in S$,
\begin{equation}\label{translated}
\begin{split}
(\mathcal{K}_{m,n}\varphi_\textit{qper})(\mathbf{x})&=e^{-i(m \mathbf{k}\cdot\mathbf{v}_1+
  n \mathbf{k}\cdot\mathbf{v}_2)}\int_{S}\mathcal{G}(\mathbf{x}-(\mathbf{x}'-m\mathbf{v}_1-n\mathbf{v}_2))\varphi(\mathbf{x}')ds(\mathbf{x}')
\\
&=e^{-i(m \mathbf{k}\cdot\mathbf{v}_1+
  n \mathbf{k}\cdot\mathbf{v}_2)}\int_{S}\mathcal{G}((\mathbf{x}+m\mathbf{v}_1+n\mathbf{v}_2)-\mathbf{x}')\varphi(\mathbf{x}')ds(\mathbf{x}').
\end{split}
\end{equation}
Clearly, the integral~\eqref{S_oper} evaluated at
$\mathbf{x}+m\mathbf{v}_1+n\mathbf{v}_2$ coincides with the last
integral in equation~\eqref{translated}, and, therefore, this last
integral is approximated closely by the equivalent-source expression
$\psi_{0,0}^{c_i,eq}(\mathbf x+m\mathbf{v}_1+n\mathbf{v}_2)$ where
$\psi_{0,0}^{c_i,eq}$ is defined in
equation~\eqref{eq:field_ac_eq_sources}. It follows that the quantity
$(\mathcal{K}_{mn}\varphi_\textit{qper})(\mathbf{x})$ can in turn be
approximated closely by
\begin{eqnarray}
\label{eq:field_ac_eq_sources-mn}
\psi_{m,n}^{c_i,eq}(\mathbf x)&:=&e^{-i(m \mathbf{k}\cdot\mathbf{v}_1+
  n \mathbf{k}\cdot\mathbf{v}_2)}\psi_{0,0}^{c_i,eq}(\mathbf
x+m\mathbf{v}_1+n\mathbf{v}_2)=\sum_{j=1}^{\frac{1}{2}M^{equiv}}e^{-i(m \mathbf{k}\cdot\mathbf{v}_1+
  n \mathbf{k}\cdot\mathbf{v}_2)}\nonumber\\
&\times&\left(\xi_{i,j}^{(m)\ell}\ G_k(\mathbf x-\mathbf x_{i,j}^{\ell}+m\mathbf{v}_1+n\mathbf{v}_2)+\xi_{i,j}^{(d)\ell}\frac{\partial G_k(\mathbf x-\mathbf x_{i,j}^{\ell}+m\mathbf{v}_1+n\mathbf{v}_2)}{\partial x_{\ell}}\right).\nonumber\\
\end{eqnarray}
Calling $\psi^{c_i,eq}(\mathbf x)$ the sum of the quantities
$\psi_{m,n}^{c_i,eq}(\mathbf x)$ over all integers $m$ and $n$, in
view of equation~\eqref{eq:field_ac_eq_sources-mn} we have
\[
\psi^{c_i,eq}(\mathbf
x):=\sum_{m,n=-\infty}^\infty\psi_{m,n}^{c_i,eq}(\mathbf
x)=\sum_{j=1}^{\frac{1}{2}M^{equiv}}\left(\xi_{i,j}^{(m)\ell}\
  G_k^\textit{qper}(\mathbf x-\mathbf
  x_{i,j}^{\ell})+\xi_{i,j}^{(d)\ell}\frac{\partial
    G_k^\textit{qper}(\mathbf x-\mathbf x_{i,j}^{\ell})}{\partial
    x_{\ell}}\right)
\]
provides  a close approximation of the quantity
\begin{equation}\label{gper_eq_srcs}
  \sum_{m,n\in\mathbb{Z}}\int_{S_{m,n}}\mathcal{G}(\mathbf{x}-\mathbf{x'})\varphi_{qper}(\mathbf{x}')ds(\mathbf{x'}),\quad \mathbf{x}\not\in\mathcal{S}_i.
\end{equation}

The approximating expression~\eqref{gper_eq_srcs} contains the
quasi-periodic Green's function $G_k^\textit{qper}$, and it at this
point that the proposed accelerated algorithm utilizes the windowed
periodic Green function: replacing $G_k^\textit{qper}$ in this
expression by its windowed approximation 
\begin{eqnarray}\label{ga_formula}
  G^a(\mathbf{x}-\mathbf{x'}) &=& \frac{1}{4\pi}\sum_{m,n\in\ZZ} \frac{e^{ik\left(\|(x_1-x_1',x_2-x_2')+m\mathbf{v}_1+n\mathbf{v}_2\|^2+(x_3-x_3')^2\right)^{1/2}}}{\big(\|(x_1-x_1',x_2-x_2')+m\mathbf{v}_1+n\mathbf{v}_2\|^2+(x_3-x_3')^2\big)^{1/2}}\nonumber\\
&\times&e^{-i(m \mathbf{k}\cdot\mathbf{v}_1+n \mathbf{k}\cdot\mathbf{v}_2)} \,\chi({\textstyle\frac{\|(x_1-x_1',x_2-x_2')+m\mathbf{v}_1+n\mathbf{v}_2\|}{a}}),
\end{eqnarray}
which, as established in Theorem~\ref{thm:super-algebraic}, gives rise
to super-algebraic convergence as $a\to +\infty$, we obtain the
corresponding super-algebraically close approximation
\begin{equation}\label{eq:all_equiv_sources}
  \psi^{c_i,eq}(\mathbf x):=\sum_{j=1}^{\frac{1}{2}M^{equiv}}\left(\xi_{i,j}^{(m)\ell}\ G^a(\mathbf x-\mathbf x_{i,j}^{\ell})+\xi_{i,j}^{(d)\ell}\frac{\partial G^a(\mathbf x-\mathbf x_{i,j}^{\ell})}{\partial x_{\ell}}\right)\  \mbox{for $(i,x)$ such that $\mathbf{x}\not\in\mathcal{S}_i$}.
\end{equation}
(Note that the $k$ dependence is explicitly displayed in the notation
$G_k$ for the free-space Green function, but, for notational
simplicity, it is suppressed in the notation $G^a$ for the windowed
periodic Green function used in, e.g.,
equation~\eqref{eq:all_equiv_sources}.) Since for a given $\ell$ the
circular regions $\Pi_i^\ell$ are not pairwise disjoint, it is
necessary, as indicated in~\cite{br-k1}, to combine equivalent source
intensities for sources supported at a given point $\mathbf x'$ that
corresponds to two different cells, say, $c_r$ and $c_s$ for which
$\mathbf x'= \mathbf x_{r,p}^{\ell} = \mathbf x_{s,q}^{\ell} $ for
some integers $p$ and $q$. We thus define the quantities
\begin{equation}
\label{eq:sum_sources}
\psi^{(*)\ell}(\mathbf x)=\sum_{\mathbf{x'}\in\Pi^{\ell}}\left(\xi_{\mathbf{x'}}^{(m)\ell}G^a(\mathbf x-\mathbf x')+
  \xi_{\mathbf{x'}}^{(d)\ell}\frac{\partial G^a(\mathbf x-\mathbf x')}{\partial x_{\ell}'}\right)
\end{equation}
where $\xi_{\mathbf{x'}}^{(m)\ell}$ and $\xi_{\mathbf{x'}}^{(d)\ell}$
denote the sum of all intensities of equivalent sources located at a
point $\mathbf{x'}\in \Pi^\ell$:
$$\xi_{\mathbf{x'}}^{(m)\ell}=\sum_{\mathbf
  x^{\ell}_{i,j}=\mathbf{x'}}\xi^{(m)\ell}_{i,j}\ \ \
\xi_{\mathbf{x'}}^{(d)\ell}=\sum_{\mathbf
  x^{\ell}_{i,j}=\mathbf{x'}}\xi^{(d)\ell}_{i,j}.$$ 

Note that, while the quantity $\psi^{(*)\ell}$ contains contributions
from cells $c_i$ for which the far-field restriction
$\mathbf{x}\not\in\mathcal{S}_i$ is not satisfied, the algorithmic
evaluation of the quantity~\eqref{eq:all_equiv_sources} does proceed
by evaluating $\psi^{(*)\ell}$ (by means of an FFT) and then
correcting for nearby contributions
$\mathbf{x}\in\mathcal{S}_i$. These two steps in the algorithm are
considered in the following  subsections.

\subsubsection{FFT evaluation of the convolutions and Correction step\label{corr}} 
As indicated above, the inaccurate quantity
$\psi^{(*)\ell}(\mathbf{x})$ (equation~\eqref{eq:sum_sources}) plays
an important role in the proposed accelerated quasi-periodic
solver. For each $\ell =1,2,3$ the proposed algorithm first evaluates
the Cartesian convolutions $\psi^{(*)\ell}(\mathbf x)$ ($\mathbf x\in
\Pi^\ell$) by means of the three-dimensional FFT algorithm.  The
proposed use of the quasi-periodic Green function, which only occurs
in the algorithm as part of the acceleration step, provides the
additional advantage that, under the strategies mentioned in
Section~\ref{sec:cost}, the Green function needs to be evaluated at a
number of the order of ${\mathcal O}(N^{4/3})$ points only---and not
for the ${\mathcal O}(N^{2})$ pairs of discretization points, where
${\mathcal O}(N)$ is the number of grid points that are used to
discretize the scatterers in the reference cell. As demonstrated in
Section~\ref{sec:numerical}, the combined windowed-Green-function
FFT-based algorithm provides a very efficient quasi-periodic
solver---at least away from Wood anomalies.

But, as indicated above, corrections are necessary to the pure
FFT-based quantity $\psi^{(*)\ell}(\mathbf{x})$: the incorrect
contributions $\mathbf{x}\in\mathcal{S}_i$ must be subtracted, and
corresponding accurate replacements need to be added. In some detail,
the quantity $\psi^{(na,eq)\ell}(\mathbf{x})$, which equals the sum of
the values at the point $\mathbf{x}\in S$ of all fields arising from
equivalent sources nonadjacent to $c_i$ can be obtained by subtracting
from $\psi^{(*)\ell}(\mathbf{x})$ the field arising at $\mathbf{x}$
from equivalent sources located within $\mathcal{S}_i$, where $i$ is
the index for which $\mathbf{x}\in c_i$.  The ``corrections''
necessary to produce $\psi^{(na,eq)\ell}(\mathbf{x})$ from
$\psi^{(*)\ell}(\mathbf{x})$ can also be evaluated efficiently, by
means of a sequence of (small) three-dimensional FFTs, since they only
involve (small) three-dimensional convolutions and {\em free-space
  Green's functions}. Once completed for $\ell =1$, 2, 3, this overall
procedure results in accurate values, on a mesh that samples the
boundaries of all cells $c_i$, of the fields arising from all true
sources contained in all cells $c_j$ not adjacent to $c_i$.

In order to obtain approximations of the nonadjacent interactions
$\psi^{(na,true)}(\mathbf{x})$ (that is, the fields generated at
$\mathbf{x}$ by the true discrete surface sources contained outside
$\mathcal{S}_i$) at surface points ${\mathbf x}\in S\cap c_i$, the
algorithm employs solutions to the Helmholtz equation within $c_i$,
with Dirichlet boundary conditions given by $\psi^{(na,eq)\ell}$,
$\ell =1$, 2, 3. These Dirichlet problems can be solved uniquely (in
view of our assumption that the wavenumber $k$ is not a resonant
frequency), and thus the good approximation properties of the
nonadjacent interactions on the boundary of each cell $c_i$ translate
into good approximations for the nonadjacent interactions on the
surface $S$. Following~\cite{br-k1}, our algorithm produces the needed
solutions of Dirichlet problems by means of approximations of the form
\begin{equation}
\label{eq:plane_wave_exp}
P(\mathbf x) = \sum_{j=1}^{n^{w}}\gamma^{j} \exp(ik\mathbf
u_j\cdot{\mathbf{x}}),
\end{equation}
valid within $c_i$ (in terms of plane wave solutions of the Helmholtz
equation), for the field $\psi^{(na,true)}$. Here $\mathbf u_{j}$ are
unit vectors that adequately sample the surface of the unit sphere,
and the coefficients $\gamma^{j}$ are obtained in such a way that the
relation $P(\mathbf x) = \psi^{(na,true)}(\mathbf x)$ is satisfied, in
the least-squares sense, for all $\mathbf x$ in an adequately chosen
collocation mesh on the cubic surface $\mathcal{S}_i$.
 
\subsubsection{Adjacent interactions\label{adj_int}} Having evaluated,
by means of FFTs and plane wave expansions, accurate approximations of
the surface values of the field $\psi^{(na,true)}(\mathbf x)$ produced
by the non-adjacent surface sources (for all discretization points
$\mathbf x\in S$), surface values of the {\em total field} are then
obtained by direct addition of necessary singular and non-singular
adjacent surface sources. Briefly, the fields that need to be added to
(the approximations just obtained for) the field
$\psi^{(na,true)}(\mathbf x)$ (for a point $\mathbf x\in S$)
include~(i) Adjacent regular sources, that is, trapezoidal-rule
contributions to the integral operator from sources lying outside the
support of the floating POU $\eta_{\mathbf{x}}$ but inside
$\mathcal{S}_i$ (none of which are included in
$\psi^{(na,true)}(\mathbf x)$), and~(ii) Adjacent singular sources,
that is, the local contributions to the integral operator considered
in stage~(a) of Section~\ref{sec:evaluation}.

\subsection{Computational cost}\label{sec:cost}
It is easy to estimate the computational cost of the proposed
windowed-Green-function/accelerated algorithm for quasi-periodic
scattering problems. Indeed, the cost of the algorithm is the same as
that of its non-periodic counterpart~\cite{br-k1} except for the fact
that, in the present case, the use of the equivalent-source
intensities requires values of the {\em quasi-periodic Green function}
$G^a$, as shown in equation~\eqref{eq:sum_sources}, instead of the use
of the free-space Green function $G_k$ in the former algorithm. (Note
that the equivalent sources themselves are obtained, even in the
present periodic context, by means of the {\em free-space Green
  function} $G_k$, as shown in
equation~\eqref{eq:field_ac_eq_sources}.)  The operation count now
proceeds simply. The algorithm~\cite{br-k1} is reported to require a
cost of $\mathcal{O}(N^{4/3}\log{N})$ operations. In addition, the
windowed-Green function accelerated algorithm requires a
precomputation of the Green function $G^a(\mathbf x)$ and its
derivatives along each coordinate direction, and at all points
$\mathbf{x}$ in the accelerator meshes $\Pi^\ell,\ell=1,2,3$. These
precomputations are performed by direct summation at a cost of
$\mathcal{O}(a^2 N^{4/3})$ operations. The overall cost of the
algorithm, including all necessary Green function evaluations, thus
amounts to the $\mathcal{O}(a^2 N^{4/3})$ precomputation cost
plus the necessary number of GMRES iterations at a cost of
$\mathcal{O}(N^{4/3})$ each.

\section{Numerical results}\label{sec:numerical}

To demonstrate the speed and accuracy of the proposed accelerated
N\"ystrom algorithm we present results of applications of this method
to problems of scattering by doubly periodic arrays of
perfectly-conducting obstacles at non-Wood configurations. For
simplicity we consider two dimensional rectangular lattices of
scatterers, that is $\mathbf{v}_1=d_1(1,0,0)$ and $\mathbf{v}_2=d_2
(0,1,0)$; results of similar quality have been produced for general
lattices.  We present two main accuracy indicators, namely certain
convergence studies on one hand, and departure from energy
conservation in the numerical solution, on the other. The latter test,
which derives from the energy conservation result satisfied by the
exact PDE solution for the perfectly conducting scatterers under
consideration---namely, that the energy flux of the incident field
must equal the sum of the energy fluxes of the reflected field and the
transmitted field---can be expressed in terms of the Rayleigh
coefficients $ B^+_{j,\ell}$ of the scattering problem:
\begin{equation}
  \sum_{(j,\ell)\in P}\gamma_{j\ell}|B^+_{j,\ell}|^2+\sum_{(j,\ell)\in P}\gamma_{j\ell}|B^-_{j,\ell}+\delta_{j,\ell}^{0,0}|^2=\gamma_{00}\,,
\end{equation}
where $P$ is the set of propagating harmonics
$P=\{(j,\ell):\|\mathbf{v}^*_{j\ell}\|<k^2\}$, where $\gamma_{j\ell}$
are defined in equation~\eqref{eq:alpha_r}. The energy defect for numerically
computed Rayleigh coefficients $\widetilde{B}_{r,s}^{\pm}$ is then
defined as
\begin{equation}\label{eps_def}
  \varepsilon=\frac{\left|\sum_{(j,\ell)\in P}\gamma_{j\ell}\big( |\widetilde{B}_{j,\ell}^+|^2+|\widetilde{B}_{j,\ell}^-+\delta_{j,\ell}^{0,0}|^2 \big) -\gamma_{0,0}\right|}{\gamma_{0,0}}
\end{equation}
(where $\delta_{j,\ell}^{0,0}$ equals 1 for $(j,\ell)=(0,0)$ and zero
otherwise). Experiments based on fully converged solutions (as
verified by means of convergence studies), suggest that the energy
defect is an excellent indicator of solution accuracy for the integral
solvers under consideration.

All of the numerical examples presented in this section concern
problems of scattering by periodic arrangements of either spherical or
bean-shaped scatterers~\cite{br-k1}, both of which have diameter equal
to $2$. In all cases the periods are given by $d_1=d_2=4$, and
plane-wave incident fields with incidence angles $\psi=\phi=0$ (that
is, normal incidence) and $\psi=\phi=\pi/3$ (oblique incidence) are
considered. For these experiments we have used the accelerator
parameters $L=3$, $M^{equiv}=4$, $n^{coll}=8$, and $n^w=4$. In all
cases the linear systems resulting from our discretization was solved
by means of the GMRES iterative solver with a relative residual
tolerance~$Tol$. The tolerance value $Tol=10^{-8}$ was used to produce
Table~\ref{high-order} while the less restrictive
``adequate-accuracy'' tolerance $10^{-4}$ was used for
Tables~\ref{errors0} and~\ref{errors1}. Table~\ref{high-order}
showcases the high-order accuracy achieved by our periodic solvers in
the case of bi-periodic arrays of spheres under normal incidence.
Tables~\ref{errors0} and \ref{errors1} present results for periodic
arrays of spheres and bean-shaped obstacles for various wavenumbers
$k$ and various values of the window-radius $a$.

The errors $\varepsilon$ presented in these tables was evaluated in
accordance with equation~\eqref{eps_def}. The error $\varepsilon_1$,
on the other hand, was calculated as the absolute error in the
Rayleigh coefficient $B^+_{0,0}$ (as estimated by comparison with a
reference solution obtained by means of a highly-refined
discretization, a large value $a$ and a sufficiently small tolerance
$Tol$). We also report numbers of iterations and computational times
required by the GMRES solvers to reach the tolerance $Tol$ in each
case. The results were obtained by means of a C++ implementation of
our solvers on a single core of a 2.67 GHz Intel Xeon CPU with 24Gb of
RAM.

\section{Conclusions}\label{sec:concl}
This paper demonstrates that the previous two-dimensional windowed
Green-function methodology~\cite{BrunoDelourme} for quasi-periodic
scattering problems can successfully be extended to the
three-dimensional context. In particular, this paper presents the
first rigorous proof of super-algebraic convergence of the windowed
Green-function method in three-dimensional space. An accelerated
windowed Green-function algorithm is presented, which possesses
excellent properties.  Comparisons, in simple examples, with one of the
most advanced techniques for evaluation of periodic Green
functions~\cite{Guerin} (which is based on a combination of
resummation and partitioning techniques) suggests that the proposed
methodology can be orders of magnitude less expensive than former
approaches.

\begin{table}
\begin{center}
\begin{tabular}{|c|c|c|c|c|c|c|c|}
\hline
Scatterer &$k$ & Unknowns & $a$ & $\varepsilon$ & $\varepsilon_1$ & Iter \\
\hline
Sphere & 1 & $6\times 64\times 64$ & 25 & 1.1 $\times$ $10^{-3}$& 1.9 $\times$ $10^{-3}$ & 11 \\
Sphere & 1 & $6\times 64\times 64$ & 50 & 1.2 $\times$ $10^{-4}$& 6.1 $\times$ $10^{-5}$ & 11  \\
Sphere & 1 & $6\times 64\times 64$ & 75 & 5.0 $\times$ $10^{-6}$& 2.1 $\times$ $10^{-6}$ & 11 \\
Sphere & 1 & $6\times 64\times 64$ & 150 & 3.8 $\times$ $10^{-7}$& 3.5 $\times$ $10^{-8}$ & 11 \\
\hline
\end{tabular}
\caption{\label{high-order} Convergence of the periodic solvers using $G^{a}$ for increasing values of the truncation radius $a$ for doubly periodic arrays of spheres under normal incidence. The reference solution corresponds to $a=400$, in which case the conservation of energy error was $\varepsilon=1.2\times 10^{-7}$.}
\end{center}
\end{table}

\begin{table}
\begin{center}
\begin{tabular}{|c|c|c|c|c|c|c|c|c|c|c|}
\hline
Scatterer &$k$ & N & $a$ & $\varepsilon$ & $\varepsilon_1$ & Iter & \multicolumn{3}{c|}{Computational Times}\\
\cline{8-10}
& & & & & & & Set-up & Time/It & Total  \\
\hline
Sphere & 0.75 & $1350$ & 20 & 5.0 $\times$ $10^{-3}$& 6.4 $\times$ $10^{-3}$& 5 & 14sec & 0.4sec & 16sec \\
Sphere & 0.75 & $1350$ & 30 & 4.7 $\times$ $10^{-4}$& 1.6 $\times$ $10^{-3}$& 5 & 29sec & 0.4sec & 31sec \\
Sphere & 0.75 & $1350$ & 40 & 2.4  $\times$ $10^{-5}$& 2.2 $\times$ $10^{-4}$ & 5 & 51sec & 0.4sec & 53sec \\
\hline
Sphere & 9 & $5766$ & 20 & 5.0 $\times$ $10^{-3}$& 3.6 $\times$ $10^{-3}$ & 13 & 14sec & 3.4sec & 57sec \\
Sphere & 9 & $5766$ & 30 & 1.1 $\times$ $10^{-3}$& 1.3 $\times$ $10^{-3}$ & 13 & 29sec & 3.4sec & 1m14sec \\
Sphere & 9 & $5766$ & 40 & 7.0 $\times$ $10^{-5}$& 2.1 $\times$ $10^{-4}$ & 13 & 51sec & 3.4sec & 1m35sec \\
\hline
\hline
Bean & 0.75 & $1350$ & 20 & 3.3 $\times$ $10^{-3}$& 5.5 $\times$ $10^{-3}$ & 10 & 14sec & 1.2sec & 26sec \\
Bean & 0.75 & $1350$ & 30 & 1.9 $\times$ $10^{-3}$& 1.4 $\times$ $10^{-3}$ & 10 & 29sec & 1.2sec & 42sec \\
Bean & 0.75 & $1350$ & 40 & 3.2 $\times$ $10^{-4}$& 3.4 $\times$ $10^{-4}$ & 10 & 51sec & 1.2sec & 1m5sec \\
\hline
Bean & 9 & $5766$ & 20 & 6.1 $\times$ $10^{-3}$& 4.0 $\times$ $10^{-3}$ & 17 & 14sec & 5.35sec & 1m45sec \\
Bean & 9 & $5766$ & 30 & 1.1 $\times$ $10^{-3}$& 9.9 $\times$ $10^{-4}$ & 17 & 29sec & 5.35sec & 2m0sec \\
Bean & 9  & $5766$ & 40 & 3.2 $\times$ $10^{-5}$& 1.7 $\times$ $10^{-4}$ & 17 & 51sec & 5.35sec & 2m30sec \\
\hline
\end{tabular}
\caption{\label{errors0} Convergence of the periodic solvers using $G^{a}$ for increasing values of the truncation radius $a$ for doubly periodic arrays of spherical and bean-shaped scatterers under normal incidence.}
\end{center}
\end{table}

\begin{table}
\begin{center}
\begin{tabular}{|c|c|c|c|c|c|c|c|c|c|c|}
\hline
Scatterer &$k$ & $N$ & $a$ & $\varepsilon$ & $\varepsilon_1$ & Iter & \multicolumn{3}{c|}{Computational Times}\\
\cline{8-10}
& & & & & & & Set-up & Time/It & Total  \\
\hline
Sphere & 9 & $5766$ & 20 & 8.0 $\times$ $10^{-3}$& 5.2 $\times$ $10^{-3}$ & 23 & 14sec & 3.4sec & 1m31sec \\
Sphere & 9 & $5766$ & 30 & 3.7 $\times$ $10^{-3}$& 8.0 $\times$ $10^{-4}$ & 22 & 29sec & 3.4sec  & 1m44sec \\
Sphere & 9 & $5766$ & 50 & 4.5 $\times$ $10^{-5}$& 1.7 $\times$ $10^{-4}$ & 22 & 1m25sec & 3.4sec & 2m40sec \\
\hline
\hline
Bean & 9 & $5766$ & 20 & 4.4 $\times$ $10^{-3}$& 7.8 $\times$ $10^{-3}$ & 21 & 14sec & 5.35sec & 2m6sec \\
Bean & 9 & $5766$ & 30 & 1.2 $\times$ $10^{-3}$& 3.1 $\times$ $10^{-3}$ & 21 & 29sec & 5.35sec & 2m23sec \\
Bean & 9  & $5766$ & 50 & 3.0 $\times$ $10^{-5}$& 2.1 $\times$ $10^{-4}$ & 21 & 1m25sec & 5.35sec & 3m17sec \\
\hline
\end{tabular}
\caption{\label{errors1} Convergence of the periodic solvers using
  $G^{a}$ for increasing values of the truncation radius $a$ for
  doubly periodic arrays of spherical and bean-shaped scatterers under
  oblique incidence $\phi=\psi=\pi/3$.}
\end{center}
\end{table}

\bigskip
\noindent
{\bfseries\large Acknowledgments.}\
The authors gratefully acknowledge support from  AFOSR and NSF under contracts FA9550-15-1-0043 and DMS-1411876 (OB);  NSF DMS-0807325 (SPS); NSF DMS-1008076 (CT); and NSF DMS-0707488 and NSF DMS-1211638 (SV).

\end{document}